\documentclass[a4j,12pt]{article}

\usepackage{amsmath}
\usepackage{amssymb}
\usepackage{amsthm}
\usepackage{graphics}
\usepackage[dvips]{graphicx}
\usepackage{amscd}

\theoremstyle{definition}

\newtheorem{example}{example}%[section]

\newtheorem{prop}[example]{Proposition}
\newtheorem{rem}[example]{Remark}

{\par\ifdim\lastskip<\bigskipamount\removelastskip\fi\smallskip%
\noindent\begin{trivlist}%
\item[\hskip\labelsep{\bf Proof\ }]}%
{\rule[-.1mm]{1.5mm}{2ex}\end{trivlist}\medskip\par}

\newtheorem{theorem}[example]{Theorem}
\theoremstyle{remark}

\def\Real{{\Bbb R}}

\def\Rc{\textrm{Ric}}
\def\rc{\textrm{Ric}}

\theoremstyle{remark}

\begin{document}
\title{Lorentz Ricci solitons on 3-dimensional Lie groups}
\author{Kensuke Onda\footnote{kensuke.onda@math.nagoya-u.ac.jp}}
\date{Graduate School of Mathematics, Nagoya University }
% Nagoya, 464-8602, Japan \\
%  '09.5.30(Stu)  }
 \maketitle
%  (PRELIMINARY VERSION 0.958)
% kensuke.onda@math.nagoya-u.ac.jp
%  Keywords : Lie group, Lorentzian left-invariant metric, Ricci soliton, Heisenberg group
 
% NOT of DISTRIBUTRION  
% Ryoichi Kobayashi,

% class: 53C22, 53C30, 53C50, 22E99

\section*{Abstract}
The three-dimensional Heisenberg group $H_3$ has three left-invariant Lorentz metrics $g_1$  , $g_2$ , and $g_3$ as in $\cite{R92}$ . 
They are not isometric each other.
In this paper, we characterize the left-invariant Lorentzian metric $g_1$ as a Lorentz Ricci soliton. 
This Ricci soliton $g_1$ is a shrinking non-gradient Ricci soliton.
Likewise we prove that the isometry group of flat Euclid plane $E(2)$ 
and the isometry group of flat Lorentz plane $E(1, 1)$ have Lorentz Ricci solitons.

\section{Introduction}
In this paper, we consider Lorentz metrics on $3$-dimensional Lie groups.
% In particular,  we take up a story about the Heisenberg group $H_3$ here.
The Heisenberg group $H_3$ has three left-invariant Lorentzian metrics $g_1$ , $g_2$ , and $g_3 \ .$ 
Rahmani $\cite{R92}$ proved that these metrics are non-isometric each other, and 
Turhan $\cite{T08} $ showed that these metrics are geodesically complete.
Rahmani $\cite{R92, RR06}$ exhibited that 
the Lie algebra of infinitesimal isometries of $(H_3, g_1)$ and $(H_3, g_2)$ are four-dimensional 
and solvable but not nilpotent, 
and the Lie algebra of infinitesimal isometries of $(H_3, g_3)$ is six-dimensional.
Moreover Rahmani $\cite{RR06}$ showed that 
the left-invariant Lorentzian metric $g_2$ has a negative constant curvature $-1/4$ , 
  $g_3$ is flat, and $g_1$ is not Einstein.
There are a difference of a Riemannian case and a Lorentzian case, and one of interesting phenomenon of a Lorentzian case. 
In this paper, we characterize the left-invariant Lorentzian metric $g_1$ as a Lorentz Ricci soliton, 
that is one of the generalization of Einstein metrics.

Let $g_0$ be a pseudo-Riemannian metric on manifold $M^n$.
If $g_0$ satisfies
$$2\rc [g_0] + L_X g_0 + \alpha g_0 = 0 \ ,$$ 
where $X$ is some vector field and $\alpha$ is some constant, 
then $(M^n, g_0, X, \alpha )$ is called a {\it Ricci soliton structure} and $g_0$ is called {\it the Ricci soliton.} 
Moreover we say that the Ricci soliton $g_0$ is a {\it gradient Ricci soliton} 
if the vector field $X$ satisfies $X=\nabla f$, where $f$ is some function, 
and the Ricci soliton $g$ is a {\it non-gradient Ricci soliton} 
if the vector field $X$ satisfies $X\ne \nabla f$ for any function $f$. 
If a constant $\alpha$ is negative, zero, or positive,  
then $g$ is called a shrinking, steady, or expanding Ricci soliton, respectively.
According to \cite{CK}, Ricci solitons have a relation with the Ricci flow.
\begin{prop}
{\it A pseudo-Riemannian metric $g_0$ is a Ricci soliton if and only if 
$g_0$ is a solution of the Ricci flow equation,}  
$$\frac{\partial }{\partial t} g(t)_{ij} = -2\Rc [g(t)]_{ij} \ ,$$
{\it satisfying 
$g(t) = c(t)  (\varphi _t)^* g_0$, where $c(t)$ is a scaling parameter, 
and $\varphi _t$ is a diffeomorphism.} 
\end{prop}
% As it is known, a metric $g_0$ is a Ricci soliton if and only if 
% $g_0$ is a solution of the Ricci flow equation,  
% $$\frac{\partial }{\partial t} g(t)_{ij} = -2\Rc [g(t)]_{ij} \ ,$$
% satisfying 
% $g(t) = c(t)  (\varphi _t)^* g_0$, where $c(t)$ is a scaling parameter, 
% and $\varphi _t$ is a diffeomorphism.
An interesting example of Ricci solitons are $(\Real ^2 , g_{st} \ , X, \alpha )$, 
where the metric $g_{st}$ is the Euclidean metric on $\Real ^2$, 
the vector field $X$ is 
$$X= -\frac{1}{2} (\alpha x + \beta y) \frac{\partial }{\partial x} +\frac{1}{2} (\beta x -\alpha y) \frac{\partial }{\partial y} \ , $$
and $\alpha $ and $\beta$ are any constants.
If $\beta =0$, then $(\Real ^2 , g_{st} \ , X, \alpha )$ are gradient Ricci solitons, named Gaussian solitons, 
and if $\beta \ne 0$, then $(\Real ^2 , g_{st} \ , X, \alpha )$ are non-gradient Ricci solitons.
% \footnote{Likewise $(\Real ^n , g_{st} \ , X, \alpha )$ has a gradient Ricci soliton structure and a non-gradient Ricci soliton structure.}  
In closed Riemannian case, Perelman $\cite{P02}$ proved that Ricci solitons become gradient Ricci solitons, 
and any steady or expanding Ricci solitons become Einstein metrics with Einstein constant zero or negative, respectively.
But in non-compact Riemannian case, it is false.
A counterexample is that 
any left-invariant Riemannian metrics on the three-dimensional Heisenberg group 
are expanding non-gradient Ricci solitons but not Einstein metrics $(\cite{BD07}, \cite{GIK06}, \cite{L07})$.

Our goal is to characterize the left-invariant Lorentzian metric $g_1$ as a Lorentz Ricci soliton, 
and to prove that the isometry group of flat Euclid plane $E(2)$ and the isometry group of flat Lorentz plane $E(1, 1)$
 have Lorentz Ricci solitons. 
This paper is organized as follows.
Section $2$ contains a brief review of the Levi-Civita connection 
of a left-invariant metric to prepare for Section $3$, $4$, and $5$.
In Section $3$, we calculate Ricci tensor of the left-invariant Lorentzian metric $g_1$ on $H_3$, 
and prove that $g_1$ satisfies Ricci soliton equation.
In Section $4$,  we calculate Ricci tensors of a left-invariant Lorentzian metric on $E(2)$, 
and prove that the metric satisfies a Ricci soliton equation.
Finally, in Section $5$ we consider a left-invariant Lorentzian metric on $E(1, 1)$, 
that is geodesically incomplete and satisfies a Ricci soliton equation.

\section{Preliminaries}
Let $G$ be a $3$-dimensional Lie group, 
and ${\mathfrak g}$ the Lie algebra of all left-invariant vector fields on $G$.
Suppose that $X$, $Y$, and $Z$ are left-invariant vector fields belonging to ${\mathfrak g}$.
In $\cite{CK}$, 
a map ${\rm ad} : {\mathfrak g} \rightarrow {\rm gl} ({\mathfrak g}) $ is defined by
$$({\rm ad} X) Y := [X, Y] \ , $$
and the adjoint ${\rm ad}' : {\mathfrak g} \rightarrow {\rm gl} ({\mathfrak g}) $ with respect to a left-invariant metric $g$ of 
the map ${\rm ad} : {\mathfrak g} \rightarrow {\rm gl} ({\mathfrak g}) $ given by
$$g\big( ({\rm ad} X)' Y, Z\big) = g\big( Y, ({\rm ad} X) Z\big) = g\big( Y, [X, Z]\big) \ .$$
As known, the Levi-Civita connection is determined as follow.
\begin{prop}[\cite{CK}]\label{levi}
{\it The Levi-Civita connection of a left-invariant metric $g$ is determined by ${\rm ad }$ and ${\rm ad}'$ via the formula}
\begin{equation}
\nabla _X Y =\frac{1}{2} \big \{   [X, Y] -({\rm ad} X)' Y - ({\rm ad} Y)' X \big \} \ .
\end{equation}
\end{prop}
The above proposition is proved for a pseudo-Riemannian metric.
We use the formula as later.

%%%%%%%%%%%%%%%%%%%%%%%%%%%%%Heisenberg group, nil

\section{The three-dimensional Heisenberg group}
% The Heisenberg group $H_3$ is isomorphic to the set of upper-triangular $3\times 3$ matrices
% \begin{eqnarray}
% H_3 \cong \Biggl\{ \left(
% \begin{array}{ccc}
% 1 & x & z \\
% 0 & 1 & y \\
% 0 & 0 & 1 \\
% \end{array}
% \right) : x, y, z\in \Real \Biggr\}
% \end{eqnarray} 
% endowed with the usual matrix multiplication.
% And $H_3$ is diffeomorphic to $\Real ^3$ under the map 
% \begin{eqnarray}
% H_3 \ni \left(
% \begin{array}{ccc}
% 1 & x & z \\
% 0 & 1 & y \\
% 0 & 0 & 1 \\
% \end{array}
% \right) \mapsto \Real ^3 \ .
% \end{eqnarray} 
% Lie algebra of $H_3$ has a basis consisting of 
% \begin{eqnarray}
% x= \left(
% \begin{array}{ccc}
% 0 & 1 & 0 \\
% 0 & 0 & 0 \\
% 0 & 0 & 0 \\
% \end{array}
% \right) \ , \ 
% y= \left(
% \begin{array}{ccc}
% 0 & 0 & 0 \\
% 0 & 0 & 1 \\
% 0 & 0 & 0 \\
% \end{array}
% \right) \ , \ 
% x= \left(
% \begin{array}{ccc}
% 0 & 0 & 1 \\
% 0 & 0 & 0 \\
% 0 & 0 & 0 \\
% \end{array}
% \right) 
% \end{eqnarray} 
% for which
% $$[x, y]=z \ , \ [y, z]=0 \ , \ [z, x]=0\ .$$
The Heisenberg group $H_3$ has three left-invariant Lorentzian metrics $g_1 \ , \ g_2 ,$  and $g_3 \ ,$
giving by  
\begin{align*}
g_1 = & -dx^2 +dy^2 +(x\ dy +dz)^2 \ , \\
g_2 = & dx^2 +dy^2 -(x\ dy +dz)^2 \ , \\
g_3 = & dx^2 +(x \ dy +dz)^2 -\big( (1-x) dy -dz\big) ^2 \ .
\end{align*}
Rahmani $\cite{RR06}$ showed that the left-invariant Lorentzian metric $g_2$ has a negative constant curvature $-1/4$ , 
  $g_3$ is flat, and  $g_1$ is not Einstein.
In this section, we characterize the left-invariant Lorentzian metric $g_1$ as a Lorentz Ricci soliton. 
Let our frame be defined by
$$F_1  = \frac{\partial }{\partial z} \ , 
F_2  = \frac{\partial }{\partial y} -x\frac{\partial }{\partial z} \ , 
F_3  = \frac{\partial }{\partial x} \ ,$$
and coframe 
$$\theta ^1  = x\ dy + dz \ , 
\theta ^2  = dy \ , 
\theta ^3  = dx  \ .$$
It is easy to check that the metric $g_1$ is represented by $g_1 =(\theta ^1)^2+(\theta ^2)^2-(\theta ^3)^2 \ $, 
and all brackets $[F_i , F_j ]$ vanish except $[F_2, F_3] = F_1$ . 
Then the Levi-Civita connection of the metric $g_1$ is given by
\begin{eqnarray}\label{levi2}
(\nabla _{F_i} F_j)=\frac{1}{2} \left(
\begin{array}{cccc}
0 & F_3 & F_2 \\
F_3 & 0 & F_1 \\
F_2 & -F_1 & 0 \\
\end{array}
\right) , 
\end{eqnarray} 
and its Ricci tensor is expressed by
$$ R_{11}  = -\frac{1}{2} \ ,     R_{22}  = \frac{1}{2} \ ,   R_{33}  = -\frac{1}{2} \ ,$$
and other components are $0$.
Since $R_{11}  = -\frac{1}{2} g_{11} $ and $R_{22}  = \frac{1}{2} g_{22}$, 
$g_1$ is not Einstein. 
But $g_1$ satisfies a Ricci soliton equation as follows.

\begin{theorem}\label{sLRS}
{\it The left-invariant Lorentzian metric}
$$g_1= -dx^2 +dy^2 +(x\ dy +dz)^2$$
{\it satisfies a Ricci soliton equation}
$$2\Rc + L_X g -3 g = 0 \ ,$$
{\it where the vector field $X$ is defined by}
$$X=(2 z + x y) F_1 +y F_2 +x F_3 \ .$$
{\it Moreover the vector field $X$ satisfies $X \ne \nabla f$ for any function $f$.
Therefore the left-invariant Lorentzian metric $g_1$ is a shrinking non-gradient Lorentz Ricci soliton.}
\end{theorem}

\proof 
% A calculation shows that the coordinates $(\nabla _i X_j)$ of 
Using $(\ref{levi2})$,  we get
\begin{eqnarray}
(\nabla _i X_j)=\left(
\begin{array}{cccc}
2 & \frac{1}{2} x & -\frac{1}{2} y \\
-\frac{1}{2} x & 1 & -\frac{1}{2} (2 z + xy) \\
\frac{1}{2} y & \frac{1}{2} (2 z + x y) & -1 \\
\end{array}
\right) . 
\end{eqnarray} 
Hence the metric $g_1$ is a shrinking Ricci soliton.
Since $$\nabla _1 X_2 -\nabla _2 X_1 =x \ne 0 \ , $$
we obtain the statement.
\endproof
 
As a consequence, Lorentz metrics on $H_3$ are negative constant curvature metric, 
flat metric, or shrinking Ricci soliton, which is one of the generalization of positive Einstein metrics.
 
 \begin{rem}
According to Theorem $\ref{sLRS}$, the left-invariant Lorentzian metric $g_1$ is a ``shrinking" Ricci soliton.
But any left-invariant Riemannian metrics on $H_3$ are ``expanding" Ricci solitons.
This phenomenon is understood as below.

We consider left-invariant Riemannian metrics $g(t)$ on $H_3$, giving by 
$g(t) =A(t)  (\theta ^1)^2 +B(t)  (\theta ^2)^2 +C(t) (\theta ^3)^2 \ $ ,
where it's frame $\{ F_i\} _{i=1}^{3}$ satisfies $[F_2, F_3] = F_1$, 
and $A(t)$, $B(t)$, and $C(t)$ are positive constants.
For any $A$, $B$, and $C$, the metric $g(t)$ becomes an expanding Ricci soliton.
The Ricci flow of $g(t)$ is given by  
\begin{equation}
    \begin{cases}
    \dfrac{d}{dt} A = -\dfrac{A^2}{BC} \ , \\
    \dfrac{d}{dt} B =\dfrac{A}{C} \ , \\
    \dfrac{d}{dt} C =\dfrac{A}{B} \ ,
     \end{cases}
\end{equation}
and the Backward Ricci flow 
$$\frac{\partial }{\partial t} g(t)_{ij} = 2\Rc [g(t)]_{ij} $$
of $g(t)$ is expressed by 
\begin{equation}
    \begin{cases}
    \dfrac{d}{dt} A = \dfrac{A^2}{BC} \ ,  \\
    \dfrac{d}{dt} B =-\dfrac{A}{C} \ , \\
    \dfrac{d}{dt} C =-\dfrac{A}{B} \ .
     \end{cases}
\end{equation}
The Ricci flow and Backward Ricci flow exist on $-\infty < t < T$ and $-T<t<+\infty $, respectively,  
% The parameter $t$ of the Ricci flow and the Backward Ricci flow exist on $(-\infty , T)$ and $(-T, +\infty )$, respectively,  
%  and a parameter $t$ of Backward Ricci flow exists on $(-T, +\infty )$, 
where $T$ is some constant depending only on $A(0)$, $B(0)$, and $C(0)$ .
In the Ricci flow case, if $t\rightarrow T$, then $ABC\rightarrow +\infty $.
This phenomenon means that the metric is expanding.
In the Backward Ricci flow case, if $t\rightarrow +\infty $, then $ABC\rightarrow 0 $.
This phenomenon means that the metric is shrinking.

On the other hand, we consider left-invariant Lorentzian metrics $g(t)$ on $H_3$, giving by 
$g(t) =A(t)  (\theta ^1)^2 +B(t)  (\theta ^2)^2 -C(t) (\theta ^3)^2 \ $,  
where it's frame $\{ F_i\} _{i=1}^{3}$ satisfies $[F_2, F_3] = F_1$, 
and $A(t)$, $B(t)$, and $C(t)$ are positive constants.
If $A=B=C=1$, the metric $g(t)$ becomes $g_1$.
For any $A$, $B$, and $C$, the metric $g$ becomes a shrinking Lorentzian Ricci soliton.
The Ricci flow equations of $g(t)$ are 
\begin{equation}
    \begin{cases}
    \dfrac{d}{dt} A = \dfrac{A^2}{BC} \ ,\\
    \dfrac{d}{dt} B =-\dfrac{A}{C} \ ,  \\ %\\[1mm] 
    \dfrac{d}{dt} C =-\dfrac{A}{B} \ .
     \end{cases}
\end{equation}
This equations are same as Backward Ricci flow equations of Riemannian metrics $g(t)$.
So in Lorentz case, a parameter $t$ exists on $(-T, +\infty )$, 
where $T$ is some constant depending only on $A(0)$, $B(0)$, and $C(0)$ , 
and if $t\rightarrow +\infty $, then $ABC\rightarrow 0$, 
and this phenomenon means that the metric is shrinking.
\end{rem}
 
%%%%%%%%%%%%%%%%%%%%%%%%%isometry of flat Euclid plane, sol
 
\section{The isometry group of flat Euclid plane}
In this section, we consider the isometry group of flat Euclid plane $E(2)$.
The isometry group of flat Euclid plane has at least two kinds of Lorentzian metrics as follow :
\begin{subequations}
    \begin{align*}
g_1 & =  (\sin y \ dx +\cos y \ dz)^2 +dy^2 -(\cos y \ dx -\sin y \ dz)^2  \ , \\
g_2 & =  (\sin y \ dx +\cos y \ dz)^2 -dy^2 +(\cos y \ dx -\sin y \ dz)^2  \ . \\
 \end{align*}
 \end{subequations}
According to Theorem $4.1$ of $\cite{G96}$, these metrics are geodesically complete.
The Lorentzian metric $g_2$ is expressed by
$$g_2 = -dy^2 +dx^2 +dz^2 \ .$$
Therefore the Lorentzian metric $g_2$ is a flat Lorentzian metric.
In this section, we characterize the left-invariant Lorentzian metric $g_1$ as a Lorentz Ricci soliton. 

Let our frame be defined by
$$F_1 =\sin y \frac{\partial }{\partial x} +\cos y \frac{\partial }{\partial z}  \ , \ 
F_2  =\frac{\partial }{\partial y} \ ,  \ 
F_3  =\cos y \frac{\partial }{\partial x} -\sin y \frac{\partial }{\partial z}  \ ,
$$
and coframe 
$$\theta ^1  = \sin y \ dx +\cos y \ dz \ , \ 
\theta ^2  = dy \ ,  \ 
\theta ^3  = \cos y \ dx -\sin y \ dz \ .
$$
It is easy to check that the metric $g_1$ is given by $g_1 =(\theta ^1)^2+(\theta ^2)^2-(\theta ^3)^2 \ $, 
and the bracket relations are $[F_1 , F_2 ] = -F_3$ ,  $[F_2, F_3] = -F_1$ , and $[F_3, F_1] =0$ . 
Then the Levi-Civita connection of the metric $g_1$ is represented by
\begin{eqnarray}\label{levi3}
(\nabla _{F_i} F_j) = \left(
\begin{array}{cccc}
0 & 0 & F_2 \\
F_3 & 0 & -F_1 \\
F_2 & 0 & 0 \\
\end{array}
\right) , 
\end{eqnarray} 
and its Ricci tensor vanish except $R_{22}  = 2$ .
The Lorentzian metric $g_1$ is not Einstein obviously, but satisfies a Ricci soliton equation as follows.

\begin{theorem}
{\it The left-invariant Lorentzian metric}
$$g_1 =  (\sin y \ dx +\cos y \ dz)^2 +dy^2 -(\cos y \ dx -\sin y \ dz)^2$$
{\it satisfies a Ricci soliton equation}
$$2\Rc + L_X g -4 g=0 \ ,$$
{\it where the vector field $X$ is defined by}
$$X =2(z \cos y + x \sin y ) F_1+ 2(x \cos y - z \sin y)F_3 \ .$$
{\it Moreover the vector field $X$ satisfies $X \ne \nabla f$ for any function $f$.
Thus the left-invariant Lorentzian metric $g_1$ is a shrinking non-gradient Lorentz Ricci soliton.}
\end{theorem}

\proof 
Using $(\ref{levi3})$,  we get 
\begin{eqnarray*}
(\nabla _i X_j)=\left(
\begin{array}{cccc}
2 & -2(-z \sin y + x \cos y) & 0 \\
2(-z \sin y + x \cos y) & 0 & 2(x \sin y + z \cos y) \\
0 & -2(x \sin y + z \cos y) & -2 \\
\end{array}
\right) . 
\end{eqnarray*} 
It follows that the metric $g_1$ is a shrinking Ricci soliton.
Since $$\nabla _1 X_2 -\nabla _2 X_1 = -4(-z \sin y + x \cos y) \ne 0 \ , $$
we obtain the statement.
\endproof

%%%%%%%%%%%%%%%%%%%%%sol, Sol, lorentz
 
\section{The isometry group of flat Lorentz plane}
In this section, we consider the isometry group of flat Euclid plane $E(1, 1)$.
The isometry group of flat Lorentz plane has at least two kinds of Lorentzian metrics as follow :
\begin{subequations}
    \begin{align*}
g_1 & = -dz^2 +(e^z dx + e^{-z} dy )^2 + (e^z dx - e^{-z} dy )^2 \ , \\
g_2 & =  dz^2 +(e^z dx + e^{-z} dy )^2 - (e^z dx - e^{-z} dy )^2 \ . \\
 \end{align*}
 \end{subequations}
According to Theorem $4.2$ of $\cite{G96}$, $g_1$ is geodesically incomplete, however $g_2$ is geodesically complete.
Since $g_2$ is the same as $\cite{N79}$ , $g_2$ is flat.
% Calculating along a definition, we can check the Lorentzian metric $g_2$ is a flat metric.
In this section, we characterize the left-invariant Lorentzian metric $g_1$ as a Lorentz Ricci soliton. 

Let our frame be defined by
$$
F_1 = \frac{1}{2} \Big( e^{-z} \frac{\partial }{\partial x} + e^z \frac{\partial }{\partial y} \Big) \ , \
F_2 = \frac{1}{2} \Big( e^{-z} \frac{\partial }{\partial x} - e^z \frac{\partial }{\partial y} \Big) \ , \
F_3 = \frac{\partial }{\partial z}
$$
and coframe 
$$\theta ^1  = e^z dx + e^{-z} dy \ , \ 
\theta ^2  = e^z dx - e^{-z} dy \ ,  \ 
\theta ^3  =  dz \ .
$$
It is easy to check that the metric $g_1$ is given by $g_1 =(\theta ^1)^2+(\theta ^2)^2-(\theta ^3)^2 \ $, 
and the bracket relations are $[F_1 , F_2 ] = 0$ ,  $[F_2, F_3] = F_1$ , and $[F_3, F_1] = -F_2$. 
Then the Levi-Civita connection of the metric $g_1$ is represented by
\begin{eqnarray}\label{levi4}
(\nabla _{F_i} F_j) = \left(
\begin{array}{cccc}
0 & F_3 & F_2 \\
F_3 & 0 & F_1 \\
0 & 0 & 0 \\
\end{array}
\right) , 
\end{eqnarray} 
and its Ricci tensor vanish except $R_{33}  = -2$ .
The Lorentzian metric $g_1$ is not Einstein obviously, however satisfies a Ricci soliton equation as follows.

\begin{theorem}
{\it The left-invariant Lorentzian metric}
$$g_1 = -dz^2 +(e^z dx + e^{-z} dy )^2 + (e^z dx - e^{-z} dy )^2$$
{\it satisfies a Ricci soliton equation}
$$2\Rc + L_X g -4 g=0 \ ,$$
{\it where the vector field $X$ is defined by}
\begin{eqnarray*}
X & = & 4 a (x e^z F_1 +x e^z F_2 -\frac{1}{2} F_3 ) +4 (1 - a) (y e^{-z} F_1 -y e^{-z} F_2 +\frac{1}{2} F_3 ) \\
& + & b (e^z F_1 +e^z F_2) +c (e^{-z} F_1 - e^{-z} F_2 ) \ ,
\end{eqnarray*} 
{\it where $a$, $b$, and $c$ are any constants.
 Moreover the vector field $X$ satisfies $X \ne \nabla f$ for any function $f$.
Thus the left-invariant Lorentzian metric $g_1$ is a shrinking non-gradient Lorentz Ricci soliton.}
\end{theorem}

\proof 
Using $(\ref{levi4})$,  we get the following : 
\begin{eqnarray*}
\nabla _1 X_1 & = & 2 \ , \\
\nabla _1 X_2 & = & \nabla _2 X_1 = 0 \ , \\
\nabla _1 X_3 & = & \nabla _3 X_1 = - 4 a x e^z +4(1-a)y e^{-z} - b e^z +c e^{-z} \ , \\
\nabla _2 X_2 & = & 2 \ , \\
\nabla _2 X_3 & = & \nabla _3 X_2 = -4 a x e^z -4(1-a)y e^{-z} - b e^z -c e^{-z} \ , \\
\nabla _3 X_3 & = & 0 \ . \\
\end{eqnarray*} 
It follows that the metric $g_1$ is a shrinking Ricci soliton.
Since $$\nabla _3 X_2 -\nabla _2 X_3 = 8 a x e^z -8(1-a)y e^{-z} + 2 b e^z -2 c e^{-z} \ne 0 \ , $$
we obtain the statement.
\endproof

\begin{rem}
The vector field $X$ is described by
\begin{eqnarray*}
X & = & 4 y e^{-z} F_1 -4 y e^{-z} F_2 +2 F_3 \\
   & & + 4 a Y_1 + b Y_2 + c Y_3 \ ,  
\end{eqnarray*} 
where $Y_1 = (x e^z -y e^{-z} ) F_1 +(x e^z + y e^{-z} ) F_2 - F_3 \ , Y_2 =e^z F_1 +e^z F_2 \ , Y_3 = e^{-z} F_1 - e^{-z} F_2 \ .$
Vector fields $\{ Y_i\}_{i=1}^3$ form the Lie algebra of infinitesimal isometries of $(E(1, 1), g_1)$, that is solvable.
\end{rem}

% as a result
% \section*{Acknowledgement}
% The author thanks Professor Ryoichi Kobayashi for his discussions.
% The author got the start of this paper from conversation with him.

%%%%%%%%%%%%%%%%%%%%%%%%%%%%%%%%%%%%%%
%%%%%%            ŽQl•¶Œ£
%%%%%%%%%%%%%%%%%%%%%%%%%%%%%%%%%%%%%%

\end{document}